\documentclass[a4paper]{article}
\usepackage{amsmath,amssymb,amsxtra}

\usepackage{mathtools}
\usepackage{mathrsfs}
\usepackage{latexsym}
\usepackage{bm}
\usepackage[amsmath,thref,thmmarks]{ntheorem}
\usepackage{etoolbox}
\usepackage{fancyhdr}
\usepackage{geometry}
\usepackage{titlesec}
\usepackage{tocloft}
\usepackage[
backend=biber,
style=alphabetic,
]{biblatex}
\usepackage{enumitem}
\usepackage{array}
\usepackage{tabularx}
\usepackage{tikz-cd}
\usepackage{graphicx}
\usepackage{authblk}

\definecolor{darkgreen}{rgb}{0,0.5,0}
\definecolor{darkblue}{rgb}{0,0,0.8}
\definecolor{darkred}{rgb}{0.8,0,0}
\definecolor{lightblue}{rgb}{0,0.6,0.8}
\usepackage[pdfencoding=auto,colorlinks,citecolor=darkgreen,linkcolor=darkblue,urlcolor=darkred]{hyperref}
\usepackage[capitalize]{cleveref}

\allowdisplaybreaks[2]
\newcommand{\epsl}{\varepsilon}
\newcommand{\fai}{\varphi}
\newcommand\wan[1]{\widetilde{#1}}

\newcommand{\cl}[1]{\overline{#1}}

\newcommand{\bZ}{\mathbb{Z}}
\newcommand{\bQ}{\mathbb{Q}}
\newcommand{\bR}{\mathbb{R}}

\newcommand{\bF}{\mathbb{F}}

\newcommand\abs[1]{\lvert #1 \rvert}

\newcommand\Abs[1]{\left\lvert #1 \right\rvert}
\newcommand\Spa[1]{\left\langle #1 \right\rangle}

\DeclareMathOperator{\im}{Im}

\DeclareMathOperator{\Di}{div}

\DeclareMathOperator{\ord}{ord}

\DeclareMathOperator{\Aut}{Aut}

\DeclareMathOperator{\Gal}{Gal}

\DeclareMathOperator{\Hom}{Hom}
\DeclareMathOperator{\coker}{coker}
\DeclareMathOperator{\Sel}{Sel}
\DeclareMathOperator{\image}{im}
\DeclareMathOperator{\ur}{ur}
\DeclareMathOperator{\loc}{loc}
\DeclareMathOperator{\G}{G}
\DeclareMathOperator{\GL}{GL}
\titlespacing*{\section}{0pt}{20pt}{15pt}
\titlespacing*{\subsection}{0pt}{15pt}{12pt}
\titlespacing*{\subsubsection}{0pt}{12pt}{6pt}
{
\theoremstyle{plain}
\theorempreskip{5pt}
\theorempostskip{5pt}
\newtheorem{theorem}{Theorem}[section]
\newtheorem{cor}[theorem]{Corollary}
\newtheorem{lemma}[theorem]{Lemma}
\newtheorem{prop}[theorem]{Proposition}
\newtheorem{conj}[theorem]{Conjecture}
}
{
\theorempreskip{5pt}
\theorempostskip{5pt}
\theorembodyfont{\normalfont}
\newtheorem{defi}{Definition}[section]

}
{
\theorempreskip{5pt}
\theorempostskip{5pt}
\theorembodyfont{\normalfont}
\theoremsymbol{\ensuremath{\square}}
\newtheorem*{Proof}{Proof.}
}
{
\theoremstyle{break}
\theorempreskip{12pt}
\theorempostskip{5pt}
\theorembodyfont{\normalfont}
\theoremsymbol{\ensuremath{\square}}

}
\setlist[enumerate]{itemsep=0pt,parsep=2pt,topsep=2pt}
\setlist[itemize]{itemsep=3pt,parsep=3pt,topsep=5pt}
\makeatletter
\renewcommand{\@biblabel}[1]{[#1]\hfill}
\makeatother
  \DeclareFontFamily{U}{wncy}{}
    \DeclareFontShape{U}{wncy}{m}{n}{<->wncyr10}{}
    \DeclareSymbolFont{mcy}{U}{wncy}{m}{n}
    \DeclareMathSymbol{\Sh}{\mathord}{mcy}{"58} 
\fontsize{12pt}{20pt}\selectfont
\title{On Conjectures of $\mu$-Invariant for Selmer Groups at $p=2$}
\author[1]{Zichao Lin}
\author[2]{Mulun Yin}
\affil[1]{The University of Massachusetts Amherst}
\affil[2]{Morningside Center of Mathematics; Academy of Mathematics and Systems Sciences, Chinese Academy of Science}
\begin{document}
\maketitle
\abstract{Let $E$ be an elliptic curve define over $\bQ$. Further assume $E$ has good ordinary reduction at $p=2$. In this article, we prove Greenberg's conjecture on the value of the algebraic Iwasawa $\mu$-invariant when $E$ has a rational isogeny of degree 2, and obtain a sharp absolute upper bound of the $2$-adic $\mu$-invariants for all such curves. We also prove the vanishing of the algebraic $\mu$-invariant of the fine Selmer group introduced by Coates and Sujatha at $p=2$.}
\tableofcontents
\setlength{\abovedisplayskip}{5pt}
\setlength{\abovedisplayshortskip}{3pt}
\setlength{\belowdisplayskip}{5pt}
\setlength{\belowdisplayshortskip}{3pt}
\setlength{\jot}{3pt}
\section{Introduction}
\large
\subsection{Main Results}
Let $E$ be an elliptic curve defined over $\bQ$. Fix a prime $p$ where $E$ has good ordinary reduction. In \cite{Ma72}, Mazur initiated the study of the growth of the $p^\infty$-Selmer group over $\bQ_\infty/\bQ$, where $\bQ_\infty$ is the unique Galois extension of $\bQ$ with $\Gal(\bQ_\infty/\bQ)\simeq \bZ_p$. Let $\Gamma=\Gal(\bQ_\infty/\bQ)$ and $\Lambda=\bZ_p[[\Gamma]]$, which is usually identified with a formal power series ring $\bZ_p[[T]]$ when a topological generator $\gamma$ of $\Gamma$ is sent to $1+T$. Denote  by $\Sel(E[p^\infty]/\bQ_\infty)$ the $p^\infty$-Selmer group over $\bQ_\infty$, and set\[X(E[p^\infty]/\bQ_\infty)=\Hom(\Sel(E[p^\infty]/\bQ_\infty), \bQ_p/\bZ_p).\] According to \cite{Ka04}, $X(E[p^\infty]/\bQ_\infty)$ is a finitely generated $\Lambda$-torsion module. Then by the fundamental theorem of finitely generated $\Lambda$-modules (see Theorem 13.12 in \cite{Wa97}), there is a $\Lambda$-homomorphism between $X(E[p^\infty]/\bQ_\infty)$ and \[(\bigoplus_{j=1}^t \Lambda/(f_i^{m_i}))\oplus (\bigoplus_{i=1}^s \Lambda/(p^{n_i}))\] with finite kernel and cokernel, with $f_i$ some irreducible polynomials in $\Lambda$. \par
Define \[\mu_{E,p}=\sum_{i=1}^s n_i\] to be the $\mu$-invariant of $E$ with respect to $p$. In \cite{Gr99}, Greenberg made the following conjecture about $\mu_{E,p}$ among the elliptic curves in the $\bQ$-isogeny class of $E$. \par
\begin{conj}
Let $E$ be an elliptic curve defined over $\bQ$, and $p$ be a prime where $E$ has good ordinary reduction. Then there exists a $\bQ$-isogenous elliptic curve $E'$ such that $\mu_{E',p}=0$. 
\end{conj}
This is conjecture 1.11 in \cite{Gr99}.\par
In this article, we will prove the above conjecture under the hypotheses that $p=2$ and $E[2]$ is reducible as a $G_\bQ$-representation over $\bF_2$. More precisely, our first main result is the following theorem.
\begin{theorem}\label{thm1:Greenberg conj}
Let $E$ be an elliptic curve defined over $\bQ$ with good ordinary reduction at $p=2$. Further assume $E[2]$ is reducible as a $G_\bQ$-representation. Then there exists a $\bQ$-isogenous elliptic curve $E'$ such that $\mu_{E',2}=0$.
\end{theorem}
This is~\cref{thm5: Greenberg conj}. \par
Using the proof of theorem \ref{thm1:Greenberg conj}, we are able to find the upper bound of $\mu_{E,2}$ when $E$ has good ordinary reduction at 2 with $E[2]$ reducible.
\begin{theorem}
Assume $E$ is defined over $\bQ$ with good ordinary reduction at $p=2$. Further assume $E[2]$ is reducible, them $\mu_{E,2}\leq 4$.
\end{theorem}
This is~\cref{thm5: mu<=4}. In fact, this bound is sharp as 195A1 in LMFDB has $\mu_{E,2}=4$. See \cite{Gr99} for details in computation.

Besides the $\mu$-invariants of the usual Selmer groups, we will also prove a similar result for the fine Selmer groups for $p=2$. Following the ideas of Coates and Sujatha in \cite{CS05}, we can prove the following.
\begin{theorem}\label{thm1: fine mu=0}
Let $E$ be an elliptic curve defined over $\bQ$. Further suppose $E[2]$ is reducible. Let $R(E[p^\infty]/\bQ_\infty)$ be the dual of the fine Selmer group of $E[p^\infty]$ over $\bQ_\infty$. Then $R(E[p^\infty]/\bQ_\infty)$ is $\Lambda$-torsion with $\mu$-invariant 0.
\end{theorem}
This is \cref{thm5: Fine Selmer 0}.
\subsection{Strategy of the proof}
The proof of Greenberg's conjecture consists of four parts. The first part is about the relation between the $\mu$-invariants of the $2^\infty$-Selmer groups and the $4$-Selmer groups. In \cite{GV00}, Greenberg and Vatsal proved that when $p\geq 3$, $\Sel(E[p^\infty]/\bQ_\infty)^\vee/(p)$ has the same $\mu$-invariant as $\Sel(E[p]/\bQ_\infty)$ under the assumption that $E[p](\bQ)=0$. We will try to generalize this result to $p=2$ while removing the assumption that $E[2](\bQ)=0$. \par
The second ingredient is about the works by Schneider, Perrin-Riou and Drinen on the difference of $\mu$-invariant between isogeneous elliptic curves. In \cite{Sc87}, Schneider proved a formula about the difference of the $\mu$-invariants for isogenous abelian varieties, which Perrin-Riou generalized to ordinary representations in \cite{PR89}. In \cite{Dr02}, Drinen provided a formula about the relation between the difference of the $\mu$-invariants and the arithmetic properties of the kernel of the isogeny, which enables us to make explicit computations. \par
The third part is about computing $\mu_{E,2}$ under the assumption that the 2-torsion of the formal group is not $G_\bQ$-invariant. In this case we will prove that such an elliptic curve has $\mu_{E,2}=0$. To obtain this result, we need to transform the group to classical settings and apply a result in \cite{Gr77} about the $\mu$-invariant of the maximal pro-2 extension unramified outside 2 and the infinite primes. \par
In order to obtain an upper bound of $\mu_{E,2}$ under the assumption that $E[2]$ is reducible, we need to understand the possible 2-adic images of $E[2^\infty]$. In \cite{RZB15}, Rouse and Zureick-Brown provided all the possible $2$-adic images of elliptic curves defined over $\bQ$. In particular, there are elliptic curves with $\bQ$-isogeny of degree 16, while none of degree 32.
For the proof of the conjecture by Coates and Sujatha. We will just follow the ideas in \cite{CS05}. In fact, the result for $p=2$ follows after we modify some lemmas in Section 3 of \cite{CS05}.
\subsection{Relations to previous works}
The original version of the Greenberg conjecture was stated in \cite{Gr99}, in which Greenberg provided a sufficient condition for $\mu_{E,p}=0$ under the assumption that $E[p]$ is reducible. In the same article, using the property of the $\bQ$-isogeny, Greenberg conjectured that $\mu_{E,2}\leq 4$ when $E$ has good ordinary reduction at $p=2$. In \cite{GV00}, Greenberg and Vatsal provided a similar criterion to $p\geq 3 $ under the assumption that $E$ has good ordinary reduction and $E[p]$ is reducible. In other words, Greenberg conjecture holds for all these elliptic curves with the prescribed $p$. In fact, they also provided that under these assumptions, the analytic $\mu$-invariant with respect to $p$ is also 0, which offers more evidence to the Mazur's Main Conjecture.\par
After the work of Greenberg and Vatsal, the only case left out under the assumption that $E$ has good ordinary reduction and $E[p]$ is reducible is the case when $E[p]$ sits inside the nonsplit exact sequence \[0\to \bF(\varphi)\to E[p]\to \bF(\psi)\to 0.\] with $\varphi(p)\neq0$ and $\phi(-1)=1$. Under these assumptions and further assuming $\varphi$ is the trivial character, in \cite{Tr05} Trifkovic found infinitely many elliptic curves with $\mu_{E,p}=0$ where $p=3,5$. On the other hand, Hachimori proved that $\mu_{E,3}=0$ is equivalent to the vanishing of the $\mu$-invariant of a Galois group with restricted ramification over $K_\infty$, where $K$ is the real cubic subextension of $\bQ(E[3])$.\par
The study of fine Selmer group was initiated by Coates and Sujatha in \cite{CS05}, which was proposed for odd primes. In this work, we show that the fine Selmer groups have vanishing $2$-adic $\mu$-invariants when $E[2]$ is reducible.

The novelty of this work is twofold: first, probably this is the first result after \cite{Gr99} about $p=2$ with the only working assumption that $E[2]$ is reducible. Hence we are able to provide more evidence to the general case of Greenberg's conjecture. In a forthcoming article \cite{LY26}, we will try to obtain numerical criteria to pin down the elliptic curves $E$ in the fixed isogeny class with $\mu_{E,2}=0$ under the additional assumption that $E[4]$ is reducible. Second, we do not use the input from the $p$-adic $L$-functions of elliptic curves. All results that we used from the analytic side of Iwasawa theory is just the $\mu$-invariant of the class groups in \cite{FW79}.
\subsection{Further developments}
At the prime $p=2$, the remaining case is that $E[2]$ is irreducible, this time $\bQ(E[2])/\bQ$ will be a $S_3$-extension. In order to prove this, we might need more results on classical Iwasawa invariants of cyclotomic $\bZ_p$-extensions over nonabelian extensions of $\bQ$. \par
For odd primes $p$, in a forthcoming article \cite{LY26b}, we will try to prove the Greenberg's conjecture for under the assumption that $E[p]$ is reducible. We may also ask the question that whether we may extend the result to the case where $E[p]$ is irreducible. \par
Also, since $\Sel(E[p^\infty]/\bQ_\infty)^\vee$ is still a finitely generated $\Lambda$-torsion module when $E$ has multiplicative reduction at $p$. It is also natural to ask whether the same result holds for primes of multiplicative reduction. This will be carried out in a future project.
\subsection{Organization of the paper}
In section \ref{2}, we will recall the definitions of the $2^\infty$-Selmer groups and the Iwasawa invariants of finitely generated torsion modules over Iwasawa algebras. In section \ref{3}, we will use Drinen's result in \cite{Dr02} to compute the difference of $\mu$-invariants between isogenous elliptic curves in the same isogeny class. Section \ref{4} is about the proof of the upper bound of the $\mu$-invariant under the assumption that $E[2]$ is reducible while $E[4]$ is irreducible. Finally, we will give the proof of our main results in section \ref{5}.
\subsection{Acknowledgements}
The authors would like to thank Robert Pollack, Tom Weston and David Zureick-Brown for helpful discussions and advice. This work is part of the first author's forthcoming Ph.D thesis. 
\section{Selmer groups and Iwasawa invariants}\label{2}
\large
Let $E/\bQ$ an elliptic curve with good ordinary reduction at $p=2$. In this section we will define the relevant Selmer groups and Iwasawa invariants and prove some basic properties.

\subsection{Selmer groups}
Let $\bQ_\infty$ be the cyclotomic $\bZ_p$-extension of $\bQ$, and let $\eta$ be a place of $\bQ_\infty$. Consider the Kummer map \[\kappa: E(\bQ_\infty)\otimes \bQ_2/\bZ_2\to H^1(\bQ_\infty, E[2^\infty]), \] and its local version \[\kappa_\eta : E((\bQ_\infty)_\eta)\otimes \bQ_2/\bZ_2\to H^1((\bQ_\infty)_\eta, E[2^\infty]).\] 
\begin{defi}
Define the 2-primary Selmer groups to be \[\Sel(E[2^\infty]/\bQ_\infty)=\ker\Bigl(H^1(\bQ_\infty, E[2^\infty])\to \prod_\eta H^1((\bQ_\infty)_\eta, E[2^\infty])/\image(\kappa_\eta) \Bigr). \]
\end{defi} 
According to \cite{Gr89} and \cite{GV00}, we can give a equivalent definition of the 2-primary Selmer group. If $l\neq2,\infty$, we let \[H_l(\bQ_\infty, E[2^\infty])=\prod_{\eta\mid l} H^1((\bQ_\infty)_\eta, E[2^\infty]).\] Let $\eta_2$ be the unique prime of $\bQ_\infty$ above $\bQ$, fix a prime $\pi$ above $\eta_2$ in $\cl{\bQ}$ and let $I_2$ be the inertia subgroup with respect to $\pi$ in $G_{\bQ_\infty}$. Let $\wan{E}$ be the reduction of $E$ modulo 2. Define \[L_2=\ker\Bigl(H^1((\bQ_\infty)_{\eta_2}, E[2^\infty])\to H^1(I_2, \wan{E}[2^\infty])\Bigr)\] and \[H_2(\bQ_\infty,E[2^\infty])=H^1((\bQ_\infty)_{\eta_2}, E[2^\infty])/L_2.\] For infinite primes, let $\bQ_n$ the subfield of $\bQ_\infty$ with $\Gal(\bQ_n/\bQ)\simeq \bZ/p^n\bZ$. Define \[P_E^\infty(\bQ_\infty)=\prod_{v\mid \infty} H^1(\bQ_{\infty,v}, E[2])/\im(\kappa_v)).\] Fix $\Sigma$ a finite set of primes containing 2, $\infty$ and all primes where $E$ has bad reduction. Denote by $\bQ_\Sigma$ the maximal extension of $\bQ$ unramified outside $\Sigma$.
\begin{prop}\label{prop2: Equivalent definition}
An equivalent defintion of the 2-primary Selmer group over $\bQ_\infty$ is \[\Sel(E[2^\infty]/\bQ_\infty)=\ker\Bigl(H^1(\bQ_\Sigma/\bQ_\infty, E[2^\infty])\to \prod_{l\in \Sigma} H_l(\bQ_\infty, E[2^\infty])\Bigr). \]
\end{prop}
\begin{Proof}
This is in section 2 of \cite{Gr89}.    
\end{Proof}
For further computation purpose, we also need to define the non-primitive 2-primary Selmer group over $\bQ_\infty$ and the $2^n$-Selmer group. Let $\Sigma_0=\Sigma-\{ 2,\infty\}$. 
\begin{defi}
Define \[S(E[2^\infty]/\bQ_{\infty})=\Sel^{\Sigma_0}(E[2^\infty]/\bQ_\infty)=\ker\Bigl( H^1(\bQ_\Sigma/\bQ_\infty, E[2^\infty])\to \prod_{l\in \Sigma-\Sigma_0} H_l(\bQ_\infty)\Bigr).\]
\end{defi}

Define \[\kappa_n: E(\bQ_\infty)/2^nE(\bQ_\infty)\to H^1(\bQ_\infty, E[2^n])\]to be the $2^n$-Kummer map, with \[\kappa_{n,\eta}: E((\bQ_\infty)_\eta)/2^nE((\bQ_\infty)_\eta)\to H^1((\bQ_\infty)_\eta, E[2^n])\] the local version.
\begin{defi}
Define the non-primitive $2^n$-Selmer group to be \[S(E[2^n]/\bQ_\infty)=\ker\Bigl(H^1(\bQ_\infty, E[2^n])\to \prod_{\eta\in \Sigma-\Sigma_0} H^1((\bQ_\infty)_\eta, E[2^n])/\image(\kappa_{n,\eta})\Bigr). \]
\end{defi}

\subsection{Iwasawa invariants}

Let $\Lambda=\bZ_p[[\Gal(\bQ_\infty/\bQ)]]\simeq \bZ_p[[T]]$, then $\Sel(E[2^\infty]/\bQ_\infty), S(E[2^\infty])/\bQ_\infty)$, and  $\Sel(E[2^n]/\bQ_\infty)$ are naturally $\Lambda$-modules. According to \cite{Ka04}, $\Sel(E[2^\infty]/\bQ_\infty)$ is a finitely generated $\Lambda$-cotorsion module. In other words, there exists a pseudo-isomorphism \[\Sel(E[2^\infty]/\bQ_\infty)^\vee\to \Bigl(\bigoplus_{i=1}^m \Lambda/(f_i)\Bigr)\bigoplus \Bigl(\bigoplus_{j=1}^n \Lambda/(2^{a_j})\Bigr)\] with finite kernel and cokernel. We define the $\mu$-invariant of the elliptic curve $E$ to be \[\mu_E=\sum_{j=1}^n a_j.\]

The reason why we introduce the non-primitive Selmer group is that it has the same $\mu$-invariant as the usual Selmer group, and is easier to work with.

\begin{prop}\label{prop2: same mu}
$S(E[2^\infty]/\bQ_{\infty})$ is $\Lambda$-cotorsion. The $\mu$-invariants of $S(E[2^\infty]/\bQ_{\infty})^\vee$ and $\Sel_2(E/\bQ_{\infty})^\vee$ are equal.
\end{prop}
\begin{Proof}
This is Corollary 2.3 in \cite{GV00}. \par
\end{Proof}

For later use, we now define the fine Selmer group for $E$ over $\bQ_\infty$. 
\begin{defi}
Define the 2-primary fine Selmer groups to be \[\Sel_0(E[2^\infty]/\bQ_\infty)=\ker(H^1(\bQ_\infty, E[2^\infty])\to \prod_\eta H^1((\bQ_\infty)_\eta, E[2^\infty]).\]
\end{defi}
It is not hard to see that our definition agrees with the one in section 3 of~\cite{CS05}, since their $S$ contains $p$, bad places of $E$ and infinite primes, and for any finite prime outside $S$ the unramified local condition agrees with strict local condition.

\subsection{Some computational results on cohomology groups}
Now we will prove some useful results on cohomology groups which will be used in chapter \ref{4}. First, let us give an alternative description of the local conditions of $S(E[2^n]/\bQ_\infty)$. Denote \[C_2=\ker\Bigl(E[2^\infty]\to \wan{E}[2^\infty]\Bigr).\] Since $E$ has good ordinary reduction at 2, $C_2$ has $\bZ_2$-corank 1. Let $D_2=E[2^\infty]/C_2$, which is a quotient module of $E[2^\infty]$ with trivial $I_2$-action. For an infinite prime $v$ of $\bQ_\infty$, denote by $C_{\infty, v}$ the divisible part of the minus part of $E[2^\infty]$ with respect to the complex conjugation induced by $v$. Let $D_{\infty,v}=E[2^\infty]/C_{\infty,v}$. Let $\Delta$ be the discriminant of $E$.   

\begin{lemma}\label[lemma]{lem2: local expression at 2}
The local condition for $S(E[2^n]/\bQ_\infty)$ at $p=2$ is $H^1(I_2, D_2[2^n])$.
\end{lemma}

\begin{Proof}
This is a generalization of Proposition 4.4 in \cite{Ma08}. For simplicity, set $F=\bQ_{\infty,2}$, $F^{\ur}=\bQ_{\infty,2}^{\ur}$. We need to show \[\ker\Bigl(g_2: H^1(F, E[2^n])\to H^1(F, E)\Bigr)\] is equivalent to \[\ker\Bigl(f_2:H^1(F, E[2^n])\to H^1(F^{\ur}, D_2[2^n])\Bigr).\] \\
As $D_2$ is divisible and $H^0(F, D_2)$ is finite, we can see that $H^1(F^{\ur}/F, D_2)=0$. Now by Proposition 4.3 in \cite{CG96}, we have \[\ker(g_2)=\ker\Bigl(H^1(F, E[p^n])\to H^1(F, D_2)\hookrightarrow H^1(F^{\ur}, D_2)\Bigr).\] On the other hand, since $H^0(F^{\ur}, D_2)=D_2$ is divisible, we conclude that \[\ker(f_2)=\ker(H^1(F, E[2^n]\to H^1(F^{\ur}, D_2[2])\hookrightarrow H^1(F^{\ur}, D_2).\] Therefore $\ker(f_2)=\ker(g_2)$. 
\end{Proof}

Let $\Delta$ be the discriminant of $E$.

\begin{lemma}\label[lemma]{lem2: local expression at infty}
If $\Delta>0$, then the local condition for $S(E[2^n]/\bQ_\infty)$ at an infinite prime $v$ is $H^1((\bQ_\infty)_v, D_{\infty,v}[2^n])$. If $\Delta<0$, then the local condition at an infinite prime $v$ is 0.
\end{lemma}

\begin{Proof}
The local condition is just the image of $E(\bR)/2^nE(\bR)$ in $H^1(\bR, E[2^n])$. Let $c$ be a complex conjugation induced by $v$, and $\Spa{P, Q}$ be a basis of $E[4]$ over $\bZ/4\bZ$ such that $c$ acts as \[\begin{pmatrix}
-1 & t \\
0 & 1
\end{pmatrix}.\] If $\Delta>0$, then $2\mid t$, hence \[H^1(\bR, E[4])=\Spa{P,2Q}/\Spa{2P}, H^1(\bR, E[2])=\Spa{2P, 2Q}.\] On the other hand, $\image(\kappa_{1,\infty})=<2P>$ and $ \image(\kappa_{2,\infty})=\Spa{P}$. Therefore, \[H^1(\bR, E[2^n])/\image(\kappa_{n,\infty})\cong \Spa{2Q}\cong H^1(\bR, D_\infty[2^n]).\] \\
If $\Delta<0$, then $2\nmid t$, hence \[H^1(\bR, E[4])=H^1(\bR, E[2])=0.\] Therefore the local condition at $v$ is 0.
\end{Proof}

\begin{prop}\label{prop2: description of E[4]}
If $\Delta>0$, then
\[S(E[2^n]/\bQ_\infty)\cong \ker\Bigl(H^1(\bQ_\Sigma/\bQ_\infty, E[2^n])\to H^1(I_2, D_2[2^n])\times \prod_{v\mid \infty} H^1(\bR, D_{\infty,v}[2^n])\Bigr).\]
If $\Delta<0$, then \[S(E[2^n]/\bQ_\infty)\cong \ker\Bigl(H^1(\bQ_\Sigma/\bQ_\infty, E[2^n])\to H^1(I_2, D_2[2^n])\Bigr).\]
\end{prop}

\begin{Proof}
This follows immediately from~\cref{lem2: local expression at 2} and~\cref{lem2: local expression at infty}.
\end{Proof}

Next we need a modification of Proposition 2.8 in \cite{GV00}, which allows us to compare the $\mu$-invariants of $\bigl(\Sel([2^{\infty}])[2^n]\bigr)^\vee$ and $\bigl(S([2^n]))^\vee$. The main difference between our result and the one in \cite{GV00} is that we need to tackle the infinite primes for $p=2$, which does not show up in the Selmer groups for $p\geq 3$. 

\begin{prop}\label[proposition]{prop2: inside outside}
The Selmer groups $S(E[2^\infty]/\bQ_{\infty})^\vee/(2^n)$ and $S(E[2^n]/\bQ_\infty)^\vee$ have the same $\mu$-invariant.
\end{prop} 

Notice that the $\Lambda/2\Lambda$-length of $S(E[2^\infty]/\bQ_\infty)^\vee$ is just $\mu_E$.

In order to prove this proposition, we need a lemma to compare the local conditions of $S(E[2^\infty])$ and $S(E[2^n])$ at infinite primes.

\begin{lemma}\label[lemma]{lem2: infinite compatible}
The natural inclusion $E[2^n]\hookrightarrow E[2^\infty]$ induces an isomorphism \[H^1(\bR, E[2^n])/\image(\kappa_{n,\infty})\cong H^1(\bR, D_\infty).\]
\end{lemma}

\begin{Proof}
By~\cref{lem2: local expression at infty}, we have \[H^1(\bR, E[2^n])/\image(\kappa_{n,\infty})\cong H^1(\bR, D_\infty[2^n]).\] But $D_\infty[2^n]\hookrightarrow D_\infty$ induces \[H^1(\bR, D_\infty[2^n])\cong H^1(\bR, D_\infty)[2^n]\cong H^1(\bR, D_\infty)\cong \bZ/2\bZ.\]
\end{Proof}

\begin{Proof}[\cref{prop2: inside outside}]
If $\Delta>0$, we consider the following commutative diagram 
\[\begin{tikzcd}
 &0 \ar{d} & \\
 &H^0(\bQ_\infty, E[2^n]) \ar{d}& \\
S(E[2^n]/\bQ_\infty) \ar{r} \ar{d}& H^1(\bQ_{\Sigma}/\bQ_\infty, E[2^n]) \ar{r} \ar{d}& H^1(I_2, D_2[2^n])\times \prod_{\eta\mid \infty} H^1(\bR, D_\infty[2^n]) \ar{d}{\fai} \\
S(E[2^\infty]/\bQ_\infty)[2^n] \ar{r} & H^1(\bQ_{\Sigma}/\bQ_\infty, E[2^\infty])[2^n] \ar{d} \ar{r} & H^1(I_2, D_2)[2^n] \times \prod_{\eta\mid \infty} H^1(\bR, D_\infty)[2^n] \\
 & 0. & 
\end{tikzcd}\]
By~\cref{lem2: infinite compatible} and the fact $H^0(I_2, D_2)=D_2$(is a divisible set), $\fai$ is an isomorphism. Now by the snake lemma, $S(E[2^n]/\bQ_\infty)$ and $S(E[2^\infty]/\bQ_\infty)[2^n]$ differ by a finite set. Therefore, their duals have the same $\mu$-invariant. \par
If $\Delta<0$, we consider instead
\[\begin{tikzcd}
 &0 \ar{d} & \\
 &H^0(\bQ_\infty, E[2^n]) \ar{d}& \\
S(E[2^n]/\bQ_\infty) \ar{r} \ar{d}& H^1(\bQ_{\Sigma}/\bQ_\infty, E[2^n]) \ar{r} \ar{d}& H^1(I_2, D_2[2^n]) \ar{d}{\fai} \\
S(E[2^\infty]/\bQ_\infty)[2^n] \ar{r} & H^1(\bQ_{\Sigma}/\bQ_\infty, E[2^\infty])[2^n] \ar{d} \ar{r} & H^1(I_2, D_2)[2^n] \\
 & 0. & 
\end{tikzcd}\]
Now $\fai$ is still an isomorphism, so as before, $S(E[2^n]/\bQ_\infty)$ and $S(E[2^\infty]/\bQ_\infty)[2^n]$ differ by a finite set.
\end{Proof}

Finally, we will use a result in \cite{Gr99} on the $\Lambda/2\Lambda$-coranks of the $H^1$ cohomology groups to prove a result about the $\mu$-factors showing up in $S(E[2^\infty]/\bQ_\infty)^\vee$.
\begin{prop}\label[proposition]{prop2: corank 1}
If $E[4]$ is reducible, then $H^1(\bQ_\Sigma/\bQ_\infty, E[2^\infty])$ has $\Lambda$-corank 1. Moreover, $H^1(\bQ_\Sigma/\bQ_\infty, E[2^\infty])/H^1(\bQ_\Sigma/\bQ_\infty, E[2^\infty])_{\Lambda- \Di}$ has $\mu$-invariant 1 if $\Delta>0$ and 0 if $\Delta<0$. In addition, the restriction map \[H^1(\bQ_\Sigma/\bQ_\infty, E[2])\to \prod_{\eta\mid \infty} H^1(\bR, E[2]). \] is surjective. Also, if $\Theta$ is a $G_\bQ$-module cyclic of order 2, then $H^1(\bQ_\Sigma/\bQ_\infty, \Theta)$ has $\Lambda/2\Lambda$-corank 1.
\end{prop}
\begin{Proof}
This is Proposition 5.8 and Lemma 5.9 of \cite{Gr99}.
\end{Proof}

\cref{prop2: corank 1} allows us to compute the $\mu$-invariants of the cohomology groups $H^1$ and the $\mu$-factors of the Selmer groups.
\begin{prop}\label[proposition]{prop2: mu of H1}
If $E[2]$ is reducible, then $H^1(\bQ_\Sigma/\bQ_\infty, E[2^n])^\vee$ has $\mu$-invariant $1+n$.
\end{prop}

\begin{Proof}
Set $B=E(\bQ_\infty)[2^\infty]$ and consider the following exact sequence \[0\to B/2^nB\to H^1(\bQ_\Sigma/\bQ_\infty, E[2^n])\to H^1(\bQ_\Sigma/\bQ_\infty, E[2^\infty])[2^n]\to 0.\]
Since $B$ is finite, $H^1(\bQ_\Sigma/\bQ_\infty, E[2^n])^\vee$ has the same $\mu$-invariant as $H^1(\bQ_\Sigma/\bQ_\infty, \\ E[2^\infty])^\vee/(2^n)$. By~\cref{prop2: corank 1}, \[H^1(\bQ_\Sigma/\bQ_\infty, E[2^\infty])^\vee/(2^n)\sim \Lambda/(2^n)\oplus\Lambda/(2),\] which has $\mu$-invariant $n+1$.
\end{Proof}

\begin{prop}\label[proposition]{prop2: mu of Sel}
If $E[2]$ is reducible, then the $\mu$-factor of $S(E[2^\infty]/\bQ_\infty)^\vee$ is $\Lambda/(2^{\mu_E})$.
\end{prop}

\begin{Proof}
Combining~\cref{prop2: corank 1} and \cref{prop2: mu of H1}, we see that the $\mu$-invariant of $S(E[2]/\bQ_\infty)^\vee$ is less than or equal to 1. Now combining~\cref{prop2: inside outside} with the structure theorem of finitely generated $\Lambda$-modules, we know that $S(E[2^\infty]/\bQ_\infty)^\vee$ has only 1 factor carrying the $\mu$-invariant.
\end{Proof}

\section{The difference of $\mu$-invariants between isogenous curves}\label{3}
In this section, we will explicitly compute the invariant $\delta$ defined in \cite{Dr02}, which allows us to compute the difference of the $\mu$-invariants between two isogenous elliptic curves.
Let $C_2$ be the kernel of the reduction map $E[2^\infty]\to \wan{E}[2^\infty]$. Following the notations in \cite{Dr02}, for a finite $G_\bQ$-submodule $\alpha$ of $E[2^\infty]$, define $F_2^+ \alpha$ to be the intersection of $C_2$ and $\alpha$. Let $v$ be the infinite prime of $\bQ$, and define the map $\phi_v: H^1(\bQ_v, \alpha)\to H^1(\bQ_v, E[2^\infty])$. Let $\ord_2$ be the normalized 2-adic valuation with $\ord_2(2)=1$. Set $\epsl_v(\alpha)=\ord_2\abs{\ker(\phi_v)}$.
\begin{defi}
Let $\alpha$ be a finite $G_\bQ$-submodule of $E[2^\infty]$, define \[\delta(\alpha)=\ord_2(\abs{F_2^+\alpha})-\ord_2(\abs{\alpha(\bQ_v)})+\epsl_v(\alpha).\]
\end{defi}

The reason that we define this $\delta$ is due to the following result of Perrin-Riou in \cite{PR89}.
\begin{theorem}\label{thm3: PR}
Let $E$ and $E'$ be two elliptic curves defined over $\bQ$ with good ordinary reduction at $p=2$. Suppose we have the following exact sequence \[0\to \alpha\to E\to E'\to 0.\] Then we have \[\mu_E-\mu_{E'}=\delta(\alpha).\]
\end{theorem}
\begin{Proof}
See Theorem 2.2 in \cite{Dr02}.
\end{Proof}
Let $\alpha$ be an order 2 $G_\bQ$-submodule of $E[2^\infty]$. Using the same notation as in \cite{Dr02}, we say $\alpha$ is ramified if $\alpha=C_2[2]$, and not ramified otherwise. For any infinite prime $v$, let $C_{\infty,v}$ denote the divisible subgroup of $E[2^\infty]$ such that a complex conjugation induced by $v$ acts as $-1$ on it. We say $\alpha$ is odd if $\alpha=C_{\infty,v}[2]$ for all the infinite primes $v$ of $\bQ_\infty$. In fact, since $\alpha$ is $G_\bQ$-invariant, $\alpha$ is odd if there exists an infinite prime $v$ of $\bQ_\infty$ such that $\alpha=C_{\infty,v}[2]$. Otherwise, we say $\alpha$ is not odd.

The main result of this section is stated as follows.
\begin{prop}\label[proposition]{prop3: delta}
Let $\alpha$ be an order 2 $G_\bQ$-submodule of $E[2^\infty]$. Then
\begin{enumerate}
\item If $\alpha$ is odd and ramified, then $\delta(\alpha)=1$.
\item If $\alpha$ is not ramified and not odd, then $\delta(\alpha)=-1$.
\item Otherwise, $\delta(\alpha)=0$.
\end{enumerate}
\end{prop}

\begin{Proof}
For the ramified part, notice that $\ord_2(\abs{F_2^+\alpha})=1$ if and only if $\alpha$ is odd, and it is 0 if and only if $\alpha$ is even.\par
For the remaining two terms in the definition of $\delta$. Let $\Delta$ be the discriminant of $E$, and $\rho: G_{\bQ_\infty}\to \Aut(E[2^\infty]\simeq GL_2(\bZ_p)$ be the representation on $E[2^\infty]$. Suppose $\Delta>0$, then for any infinite prime $v$ of $\bQ_\infty$, the induced complex conjugation $c_v$ has image $\begin{pmatrix}
-1 & t \\
0 & 1
\end{pmatrix}$ under $\rho$ with respect to a fixed basis, say $\{P,Q\}$, with $2\mid t$. Then $c_v$ acts on $E[2]$ trivially, hence $\alpha\subset E(\bQ_v)$, and $\ord_2(\abs{\alpha(\bQ_v)})=1$. Now $C_2=\Spa{P}$, \[H^1(\bQ_v,\alpha)\simeq \alpha,\] and \[H^1(\bQ_v,E[2^\infty])\simeq \frac{P\otimes \bQ_2/\bZ_2\oplus Q/2}{P\otimes \bQ_2/\bZ_2},\] hence $\epsl_v(\alpha)=1$ if and only if $\alpha=\Spa{P/2}$ if and only if $\alpha$ is odd. \par
In the case $\Delta<0$, $c_v$ has image $\begin{pmatrix}
-1 & t \\
0 & 1
\end{pmatrix}$ under $\rho$ with $2\nmid t$. Hence on the mod-2 level, the only $G_\bQ$-invariant submodule is $\Spa{P/2}$, which implies that still $\ord_2(\abs{\alpha(\bQ_v)})=1$. Since \[H^1(\bQ_v,\alpha)\simeq \alpha\] and \[H^1(\bQ_v,E[2^\infty])=0,\] we have $\epsl_v(\alpha)=1$.
\end{Proof}

\section{Upper bound of the $\mu$-invariant}\label{4}
In this section, we will try to prove the following result.
\begin{prop}\label[proposition]{prop4: nonsplit mu=0}
Let $E$ be an elliptic curve defined over $\bQ$ with good ordinary reduction at $p=2$. Suppose $E[2]$ is reducible as a $G_\bQ$-module and $C_2[2]$ is not $G_\bQ$-invariant, then $\mu_E=0$.
\end{prop}

For the setup, let $K=\bQ(E[2])$. Since $E[2]$ is reducible and $C_2[2]$ is not $G_\bQ$-invariant, $K$ is a quadratic number field. Let $\Delta=\Gal(K/\bQ)$ with $\tau$ the nontrivial element, $\{P,Q\}$ be a basis of $E[2]$ with $P\in E(\bQ)$, and $\rho_2: \Delta\to \Aut(E[2])\simeq \GL_2(\bF_2)$ with respect to the basis $\{P,Q\}$. Then $\rho_2(\tau)=\begin{pmatrix}
1 & 1 \\
0 & 1
\end{pmatrix}$.\par

We will separate the proof of~\cref{prop4: nonsplit mu=0} into three cases.

\begin{prop}\label[proposition]{prop4: K Im mu=0}
Assume all the assumptions in \cref{prop4: nonsplit mu=0}. If $K$ is imaginary, then $\mu_E=0$.
\end{prop}

\begin{prop}\label[proposition]{prop4: K Re odd sub mu=0}
Assume all the assumptions in \cref{prop4: nonsplit mu=0}. If $K$ is real and $\Spa{P}$ is odd, then $\mu_E=0$.
\end{prop}

The proofs of these two parts rely on a result by Greenberg in \cite{Gr99}, which offers a sufficient condition for $\mu_E=0$.
\begin{theorem}\label{thm4: Greenberg}
Assume $E$ is an elliptic curve defined over $\bQ$ with good ordinary reduction at $p=2$. Suppose $\Phi\subset E(\bQ)[2]$ is odd but not ramified, or not odd but ramified, then $\mu_E=0$. 
\end{theorem}
\begin{Proof}
This is Proposition 5.14 of \cite{Gr99}.
\end{Proof}

\begin{Proof}[\cref{prop4: K Im mu=0}]
Suppose the odd part of $E[2]$ is generated by $R$. Since $\tau$ acts on the odd part of $E[2^\infty]$ by $-1$, $\tau(R)=-R=R$ as $R$ has order 2. Therefore, according to the image of $\rho_2$, we have $R=P$. Since $\Spa{P}$ is odd but not ramified, then according to \cref{thm4: Greenberg}, we have $\mu_E=0$.
\end{Proof}

\begin{Proof}[\cref{prop4: K Re odd sub mu=0}]
By the assumption, $\Spa{P}\subset E(\bQ)[2]$ is odd but not ramified. Then by \cref{thm4: Greenberg}, $\mu_E=0$ 
\end{Proof}

Now we will focus on the case where $K$ is quadratic and $\Spa{P}$ is not odd. Let $v_1$ and $v_2$ be the infinite primes of $K$. Since $\Spa{P}$ is not odd, by interchanging the two primes if necessary, we may assume $C_{\infty, v_1}[2]=\Spa{Q}$, then $C_{\infty,v_2}[2]=C_{\infty,\tau v_1}[2]=\Spa{\tau(Q)}=\Spa{P+Q}$. \par
Next we claim that 2 splits in $K$. 
\begin{lemma}\label{lem4: 2 splits in K}
Assume all the assumptions listed above. Then 2 splits in $K$.
\end{lemma}
\begin{Proof}
Suppose not, then 2 is either ramified or inert in $K$. In either case, $\tau$ is in the decomposition field of $\bQ_2$. Since $C_2[2]$ is $G_{\bQ_p}$-invariant, according to the action of $\tau$ on $E[2]$, we must have $C_2[2]=\Spa{P}$, contradicting our assumption.
\end{Proof}

Since 2 splits in $K$, we assume $2O_K=\pi\bar{\pi}$ with $\pi$ the primes we selected in order to define $S(E[2]/\bQ_\infty)$. \par
Let $K_\infty=K\bQ_\infty$, $M_\infty$ to be the maximal abelian pro-2 extension of $K_\infty$ unramified outside 2 and infinite primes. In order to prove the remaining case of \cref{prop4: nonsplit mu=0}, we need a result on the structure of $\Gal(M_\infty/K_\infty)$.

\begin{theorem}\label{thm4: Greenberg Equivariant}
Assume $K\neq \bQ(\sqrt{2})$ a real quadratic field. Let $\Delta=\Gal(K/\bQ)$. Then \[\Gal(M_\infty/K_\infty)\sim \Lambda/2\Lambda[\Delta]\] as a $\Lambda[\Delta]$-module.
\end{theorem}
\begin{Proof}
This is Proposition 8 of \cite{Gr77}. The reason we need to exclude $\bQ(\sqrt{2})$ is that $\bQ(\sqrt{2})$ is the first layer of the cyclotomic tower $\bQ_\infty/\bQ$. 
\end{Proof}
In order to apply \cref{thm4: Greenberg} to our case, we need to show that our $K$ is never $\bQ(\sqrt{2})$.
\begin{lemma}\label[lemma]{lem4: not Q(sqrt2)}
Assume $E$ is an elliptic curved defined over $\bQ$ with good ordinary reduction at $p=2$. Suppose $K:=\bQ(E[2])$ is a real quadratic field, then $K\neq \bQ(\sqrt{2})$.
\end{lemma}
\begin{Proof}
In fact, we will prove $K=\bQ(\sqrt{m})$ with $2\nmid m$. Let $y^2+a_1xy+a_3y=f(x)$ be a minimal Weierstrass form of our $E$. Using a change of variable if necessary, we have $E': y^2=g(x)$ with $g(x)\in \bQ[x]$ a cubic monic polynomial. According to chapter 3.1 of \cite{Si09}, there exists $u\in \bQ^\times$ such that $\Delta_E=u^{12}\Delta_{E'}$ and $\bQ(E[2])=\bQ(E'[2])$. Hence we just need to show that $\bQ(E'[2])\neq \bQ(\sqrt{2})$. \par
Since $\bQ(E'[2])=\bQ(E[2])$ is quadratic, then $\Delta_{E'}$ is not a square in $\bQ^\times$, as otherwise $\bQ(E'[2])=\bQ$, which implies $\bQ(\sqrt{\Delta_{E'}})$ is a quadratic field. Since $\bQ(\sqrt{\Delta_{E'}})\subset \bQ(E'[2])$, we have $\bQ(\sqrt{\Delta_E})=\bQ(\sqrt{\Delta_{E'}})=\bQ(E'[2])$. Since $E$ has good ordinary reduction at $p=2$, by the Ogg's formula about the discriminant and conductor (see 4.11.1 of \cite{Si94}), $2\nmid \Delta_E$. Therefore, $\bQ(E'[2])=\bQ(\sqrt{\Delta_E})$ is of the form $\bQ(\sqrt{m})$ with $2\nmid m$.
\end{Proof}

Define \[ \loc_2: H^1(\bQ_\Sigma/\bQ_\infty, E[2])\to H^1(I_2, E[2])\] to be the localization map at $p=2$ with $I_2$ the inertia group of a prime above $\pi$. For $w\mid \infty$ in $\bQ_\infty$, define \[ \loc_w: H^1(\bQ_\Sigma/\bQ_\infty, E[2])\to H^1(\bQ_{\infty,w}, E[2])\] to be the localization map at $w$. Now we are ready to prove the remaining case of \cref{prop4: nonsplit mu=0}.
\begin{prop}\label{prop4: K real not odd sub mu=0}
Assume $E$ is an elliptic curve defined over $\bQ$ with good ordinary reduction at $p=2$. Suppose $\Phi\subset E(\bQ)[2]$ is neither odd nor ramified and $\bQ(E[2])$ is real and quadratic. Then $\mu_E=0$.  
\end{prop}
\begin{Proof}
By \cref{prop2: inside outside} and \cref{prop2: mu of Sel}. It suffices to show $S(E[2]/\bQ_\infty)^\vee$ has $\mu$-invariant 0. By definition, we have \[S(E[2]/\bQ_\infty)=\loc_2^{-1}(H^1(I_2, C_2[2]))\bigcap (\cap_{w\mid \infty} \loc_w^{-1}(H^1(\bQ_{\infty,w}, C_{\infty,w}[2]).\] We will analyze the two parts separately. \par
For the preimage of $\loc_2$, consider the following commutative diagram 
\[\begin{tikzcd}
H^1(\bQ_\Sigma/\bQ_\infty, E[2]) \ar{d}{\loc_2} \ar{r}{f} & \Hom_\Delta(\Gal(M_\infty/K_\infty), E[2]) \ar{d}{\loc_\pi\times \loc_{\bar{\pi}}} \\
H^1(I_2,E[2]) \ar{r}{g} & H^1(I_2, E[2])\times H^1(\tau (I_2), E[2])
\end{tikzcd}\] with $f$ the restriction map, $\loc_\pi$ and $\loc_{\bar{\pi}}$ the corresponding localization maps. $g$ is the map $\sigma\to (\sigma, \tau^*(\sigma))$ with $\tau^*$ the conjugation map induced by $\tau$. According to theorem 25 in \cite{Iw73}, $I_2$ has $\mu$-invariant 1, then by \cref{thm4: Greenberg} and \cref{lem4: not Q(sqrt2)}, we may identify $I_2$ with $\Lambda/2\Lambda$ as a submodule of $\Gal(M_\infty/K_\infty)\sim \Lambda/2\Lambda[\Delta]$. According to inflation-restriction sequence, we have \[\ker(f)\simeq H^1(K_\infty/\bQ_\infty, E[2])\] and \[\coker(f)\hookrightarrow H^2(K_\infty/\bQ_\infty, E[2]).\] Hence both $\ker(f)$ and $\coker(f)$ are finite. Therefore, \[\loc_2^{-1}(H^1(I_2, C_2[2]))\sim \loc_\pi^{-1}(H^1(\Lambda/2\Lambda, \Spa{Q}))\simeq \Hom(\Lambda/2\Lambda, E[2])=: N.\] \par
Next we consider the preimage of $\prod_{w\mid\infty} \loc_w$. Similar to the previous case, we consider the following commutative diagram \[\begin{tikzcd}
H^1(\bQ_\Sigma/\bQ_\infty, E[2]) \ar{d}{\prod \loc_w"} \ar{r}{f} & \Hom_\Delta(\Gal(M_\infty/K_\infty), E[2]) \ar{d} {\prod_w\loc_w\times \prod_w \loc_{\tau w}} \\
\prod_{w\mid \infty} H^1(\bQ_{\infty,w},E[2]) \ar{r}{h} & \prod_w H^1(\bQ_{\infty,w}, E[2])\times \prod_{\tau w} H^1(\bQ_{\infty,w}, E[2])
\end{tikzcd} \] with $f$ the restriction map and $h$ the map $\sigma_w\mapsto ((\sigma_w), \tau^*(\sigma_w))$. Let $M'_\infty$ be the maximal real subfield of $M_\infty$, then $\Gal(M_\infty/M'_\infty)\sim (\Lambda/2\Lambda)(1+\tau)$. Since still $\ker(f)$ and $\coker(f)$ are finite, we have \[\prod_w\loc_w(N)\sim \prod_{w\mid \infty} (H^1(\bQ_{\infty,w}, \Spa{P}).\] Since $C_{\infty,w}[2]$ is either $\Spa{P+Q}$ or $\Spa{Q}$, we obtain that for any $w\mid \infty$ in $\bQ_\infty$ \[\loc_w(N)\bigcap H^1(\bQ_{\infty,w}, C_{\infty,w}[2]).\] \par
However, according to \cref{prop2: corank 1}, $\prod_w \loc_w$ is surjective and $(\prod_w H^1(\bQ_{\infty,w}, E[2]))^\vee$ has $\mu$-invariant 1. Thus according to \cref{prop2: mu of H1}, we have \[N\cap (\cap_{w\mid \infty} \loc_w^{-1}(H^1(\bQ_{\infty,w}, C_{\infty,w}[2])))\sim 0.\] Therefore, \[S(E[2]/\bQ_\infty)^\vee\sim 0.\]
\end{Proof}

Combining \cref{prop4: nonsplit mu=0} with \cref{prop3: delta}, we can obtain an immediate corollary for the case where $E[2]$ is reducible while $E[4]$ is irreducible.
\begin{cor}\label[corollary]{cor4: mu<1}
Assume $E$ is an elliptic curve defined over $\bQ$ with good ordinary reduction at $p=2$. Suppose $E[2]$ is reducible while $E[4]$ is irreducible. Then $\mu_E\leq 1$. $\mu_E=1$ if and only if $E[2]$ contains a $G_\bQ$-invariant odd and ramified submodule $\Phi$.
\end{cor}
\begin{Proof}
Let $\Psi$ be the largest $\G_\bQ$-invariant odd and ramified submodule of $E[p^\infty]$. Since $E[4]$ is irreducible, then $\abs{\Psi}\leq 2$. Let $E'$ be the elliptic curve in the same isogeny class with $E$ and related to $E$ by the exact sequence \[0\to \Psi\to E\to E'\to 0.\] Then by \cref{prop3: delta}, we have $\mu_E-\mu_{E'}\leq 1$ and $\mu_E-\mu_{E'}=1$ if and only if $\Psi$ has order 2. \par
Let $C_{2,E}$ and $C_{2,E'}$ be the odd part of $E[2^\infty]$ and $E'[2^\infty]$ respectively. We claim that $C_{2,E'}=C_{2,E}/\Psi$. In fact, since $I_p$ acts on $C_{2,E}/\Psi$ by the cyclotomic character, then according to the classification of \cite{Se72}, $C_{2,E}/\Psi=C_{2,E'}$. Since $C_{2,E}[4]$ is not $G_\bQ$-invariant, neither is $C_{2,E'}$. Therefore, according to \cref{prop4: nonsplit mu=0}, $\mu_{E'}=0$. Combining this with the last paragraph, we have $\mu_E\leq 1$ and $\mu_E=1$ if and only if $\Psi$ has order 2.
\end{Proof}
\section{Proof of the main theorems}\label{5}
\subsection{Proof of Greenberg's conjecture}
Now we will combine the results in the previous sections to prove the two main theorems of this article.

We first recall the following result from Rouse and Zureick-Brown.
\begin{theorem}\label{thm5: 16-isogeny}
Assume $E$ is defined over $\bQ$ and has good reduction at $p=2$. Further assume $E[2]$ is reducible, then it admits at most a 16-isogeny.
\end{theorem}

\begin{Proof}
This is the main result in \cite{RZB15}.
\end{Proof}

The first main result is Greenberg's conjecture when $E[2]$ is reducible and $2$ is a good ordinary prime.

\begin{theorem}\label{thm5: Greenberg conj}
Assume $E$ is defined over $\bQ$ with good reduction at $p=2$. Further assume $E[2]$ is reducible. Then there exists a $\bQ$-isogenous elliptic curve $E'$ such that $\mu_{E'}=0$.
\end{theorem}
\begin{Proof}
Following the notations in \cref{cor4: mu<1}, let $C_{2,E}$ be the odd part of $E[2^\infty]$.
Let $\Phi$ be the largest ramified $G_\bQ$-invariant submodule of $E[2^\infty]$. According to \cref{thm5: 16-isogeny}, such $\Phi$ exists with $\Abs{\Phi}\leq 16$. Say $\Phi$ has order $2^m$. Let $E'$ the elliptic curve in the same isogeny class with $E$, related to $E$ by the exact sequence \[0\to \Phi \to E\to E'\to 0.\] Since $C_{2,E'}[2]=C_{2,E}[2^{m+1}]/\Phi$, $C_{2,E'}[2]$ is not $G_\bQ$-invariant, otherwise $C_{2,E}[2^{m+1}]$ is $G_\bQ$-invariant, contradicting the assumption that $\Phi$ is the largest ramified invariant submodule. Therefore, by \cref{prop4: nonsplit mu=0}, we conclude that $\mu_{E'}=0$.
\end{Proof}
The second main result is a universal upper bound for the 2-adic $\mu$-invariant under the assumption that $E$ has good ordinary reduction at $p=2$ with $E[2]$ reducible.
\begin{theorem}\label{thm5: mu<=4}
Assume $E$ is defined over $\bQ$ and has good ordinary reduction at $p=2$. Further assume $E[2]$ is reducible, then $\mu_E\leq 4$.
\end{theorem}

In order to prove this theorem. We need a generalization of \cref{cor4: mu<1}.

\begin{prop}\label[proposition]{prop5: mu_E<=n}
Assume $E$ is an elliptic curve defined over $\bQ$ with good ordinary reduction at $p=2$. Further assume $E[2^n]$ is reducible while $E[2^{n+1}]$ is irreducible for $n\geq 1$. Then $\mu_E\leq n$.
\end{prop}

\begin{Proof}
Let $\Phi$ be the largest possible ramified $G_\bQ$-invariant cyclic submodule of $E[2^n]$. Then by assumption, $\Phi$ has order $\leq 2^n$. Let $E'$ be the elliptic curve in the same isogeny class as $E$ with \[ 0\to \Phi \to E\to E'\to 0.\] Then according to \cref{prop3: delta}, we have $0\leq \mu_E-\mu_{E'}\leq n$. Say $\Phi$ has order $2^m$, then by the same argument as in \cref{thm5: Greenberg conj}, $C_{2,E'}[2]=C_{2,E}[2^{m+1}]/\Phi$ is not $G_\bQ$-invariant. Hence by \cref{prop4: nonsplit mu=0}, $\mu_{E'}=0$. Therefore, we conclude that $\mu_E\leq n$.
\end{Proof}

\begin{Proof}[\cref{thm5: mu<=4}]
According to \cref{thm5: 16-isogeny}, if $E[2^n]$ is reducible while $E[2^{n+1}]$ is irreducible, then $n\leq 4$. By \cref{prop5: mu_E<=n}, we have $\mu_E\leq 4$.
\end{Proof}

Finally we will give a necessary and sufficient condition for $\mu_E=n$.
\begin{theorem}\label{thm5: mu=m}
Assume $E$ is defined over $\bQ$ and has good reduction at $p=2$. Further assume $E[2]$ is reducible, then $\mu_E=n$ if and only if the largest $G_\bQ$-invariant ramified and odd subgroup $\Phi$ has order $2^n$.
\end{theorem}

\begin{Proof}
Suppose the largest $G_\bQ$-invariant ramified and odd subgroup $\Phi$ has order $2^n$. Let $E'$ be the elliptic curve with the exact sequence \[0\to \Phi \to E\to E'\to 0.\] Then by \cref{prop3: delta}, we have $\mu_E-\mu_{E'}=n$. Since $C_{2,E'}[2]=C_{2,E}[2^{n+1}]/\Phi$ is not $G_\bQ$-invariant, we have $\mu_{E'}=0$. Therefore $\mu_E=n$.
For the other direction, assume the largest $G_\bQ$-invariant ramified and odd subgroup $\Phi$ has order $2^m$. Still consider the exact sequence \[0\to \Phi \to E\to E'\to 0.\] Then still by \cref{prop3: delta}, we have $\mu_E-\mu_{E'}=m$, so $\mu_{E'}=n-m$. However, on the other hand, according to our assumption, $C_{2,E'}[2]=C_{2,E}[2^{n+1}]/\Phi$ is not $G_\bQ$-invariant. Thus by \cref{prop4: nonsplit mu=0}, we have $\mu_{E'}=0$, which implies that $n=m$.
\end{Proof}

\subsection{Proof of the Coates-Sujatha Conjecture}
Let $R(E[p^\infty]/\bQ_\infty)=\Sel_0(E[p^\infty]/\bQ_\infty)^\vee$ be the dual of the fine Selmer group of $E[p^\infty]$ over $\bQ_\infty$. In this section, we will prove the following result.
\begin{theorem}\label{thm5: Fine Selmer 0}
Assume $E$ is an elliptic curve defined over $\bQ$. Further assume $E[2]$ is reducible. Then $R(E[p^\infty]/\bQ_\infty)$ is a finitely generated $\Lambda$-torsion module with $\mu$-invariant 0.
\end{theorem}
\begin{Proof}
In fact this result is exactly the version for $p=2$ of Corollary 3.6 in \cite{CS05}. Let $K(\bQ(E[p^\infty]))$ denote the maximal unramified abelian $2$-extension of $\bQ(E[p^\infty])$ in which every prime above 2 splits completely. Then still by Lemma 3.8 of \cite{CS05}, we have \[R(E[p^\infty]/\bQ(E[p^\infty]))=\Hom(\Gal(K(\bQ(E[p^\infty]))/\bQ(E[p^\infty])), E[p^\infty]),\] as $\bQ(E[p^\infty])$ is complex due to the Weil pairing. Since $E[2]$ is reducible,  $\bQ(E[2])$ is either $\bQ$ or a quadratic extension. In either case, according to the main theorem of \cite{FW79}, the Iwasawa $\mu=0$ conjecture holds for $\bQ(E[2])$. The rest follows from the proof of Theorem 3.4 of \cite{CS05}.
\end{Proof}
\clearpage
\addcontentsline{toc}{section}{References}
\makeatletter
\interlinepenalty=10000
\printbibliography
\makeatother
\end{document}